\documentclass{article}
\usepackage{amssymb}
\usepackage{amsmath}

\title{On the heights of algebraic points on curves over number fields}
\author{Christophe Soul\'e}
\date{\it \small CNRS and Institut des Hautes \'Etudes Scientifiques, 35 route de Chartres, \\ F-91440 Bures-sur-Yvette, France \\ soule@ihes.fr}

\begin{document}

\maketitle

Let $X$ be a semi-stable regular curve over the spectrum $S$ of the integers in a number field $F$, and $\bar L = (L,h)$ an hermitian line bundle on $X$, i.e. $L$ is an algebraic line bundle on $X$ and $h$ is a smooth hermitian metric (invariant by complex conjugation) on the restriction of $L$ to the set $X({\mathbb C})$ of complex points of $X$. In this paper we are interested in the height $h_{\bar L} (D)$ of irreducible divisors $D$ on $X$ which are flat over $S$, i.e. the arithmetic degree of the restriction of $\bar L$ to $D$.

\smallskip

First we assume that the degree $\deg (L)$ of $L$ on the generic fiber $X_F$ is positive and we denote by $\bar L \cdot \bar L \in {\mathbb R}$ the self-intersection of the first arithmetic Chern class of $\bar L$. Define
$$
e(\bar L , d) = \underset{\deg (D) = d}{\inf} \ \frac{h_{\bar L} (D)}{d} \, .
$$
Our first result (Theorem 2) is that
$$
\lim_d \inf e(\bar L , d) \geq \frac{\bar L \cdot  \bar L}{2\deg (L)} \, .
$$
This is a generalization of an inequality of S.~Zhang (\cite{Z} , Th.~6.3).

\smallskip

Next, when $X_F$ has genus at least two and $\bar\omega$ denotes the relative dualizing sheaf of $X$ over $S$ with its Arakelov metric \cite{A}, we obtain in Theorem~3 explicit lower bounds for $e(\bar\omega , d)$.

\smallskip

We prove also some upper bounds. Assume that $\deg (L) > 0$ and that $\deg (L_{\vert E}) \geq 0$ for every vertical irreducible divisor $E$
on $X$. For any integer $d_0 > 0$ we define
$$
e' (\bar L , d_0) = \underset{D_0}{\rm sup} \ \underset{D \pitchfork D_0}{\rm inf} \ \frac{h_{\bar L} (D)}{\deg (D)} \, ,
$$
where $D_0$ runs over all irreducible horizontal divisors of degree $d_0$, and $D$ runs over all such divisors which meet $D_0$ properly. We prove in Theorem~4 that
$$
\lim_{d_0} \sup e' (\bar L , d_0) \leq \frac{\bar L \cdot \bar L}{2 \deg (L)} \, ,
$$
and, when $X_F$ has genus at least two, we give in Theorem~5 explicit upper bounds for $e' (\bar\omega , d_0)$.

\smallskip

The main tool in the proof of these inequalities is the lower bounds for successive minima of the lattice $H^1 (X , M^{-1})$ with its $L^2$-metric  which we obtained in previous papers \cite{Vanishing} \cite{Secant} \cite{Stable}. From these lower bounds we deduce upper bounds for the successive minima of $H^0 (X , M \otimes \omega)$ by using a transference theorem relating the successive minima of a lattice with those of its dual (Theorem~1).

\section{ Duality and successive minima : }\label{sec1 }

\subsection{ \ }\label{ssec1.1}

Let $F$ be a number field, ${\mathcal O}_F$ its ring of integers and $S = {\rm Spec} ({\mathcal O}_F)$. Consider an hermitian vector bundle $\bar E = (E,h)$ on $S$, i.e. $E$ is a finitely generated projective ${\mathcal O}_F$-module and, for every complex embedding $\sigma : F \to {\mathbb C}$, the corresponding extension $E_{\sigma} = E \underset{{\mathcal O}_F}{\otimes} {\mathbb C}$ of $E$ from ${\mathcal O}_F$ to ${\mathbb C}$ is equipped with an hermitian scalar product $h_{\sigma}$. Furthermore, we assume that $h = (h_{\sigma})$ is invariant under complex conjugation.

\smallskip

We are interested in (the logarithm of) the successive minima of $\bar E$. Namely, for any positive integer $k \leq N$, where $N$ is the rank of $E$, we let $\mu_k (\bar E)$ be the infimum of the set of real numbers $\mu$ such that there exist $k$ vectors $e_1 , \ldots , e_k$ in $E$ which are linearly independent in $E \otimes F$ and such that, for every complex embedding $\sigma : F \to {\mathbb C}$ and for all $i = 1 , \ldots , k$,
$$
\Vert e_i \Vert_{\sigma} \leq \exp (\mu) \, ,
$$
where $\Vert \cdot \Vert_{\sigma}$ is the norm defined by $h_{\sigma}$. We shall compare the successive minima of $\bar E$ with those of its dual $\bar E^*$.

\smallskip

Let $r_1$ (resp. $r_2$) be the number of real (resp. complex) places of $F$,  $r = [F:{\mathbb Q}]$ the degree of $F$ over ${\mathbb Q}$, and $\Delta_F$ its absolute discriminant.
We define
\begin{equation}
C(N,F) = \frac{1}{r} \log \vert \Delta_F \vert + \frac{3}{2} \log (N) + \frac{5}{2} \log (r) - \frac{r_2}{r} \log (\pi) \, .
\end{equation}

\bigskip

\noindent {\bf Theorem 1.} {\it For every $k \leq N$ the following inequalities hold:}
$$
0 \leq \mu_k (\bar E) + \mu_{N+1-k} (\bar E^*) \leq C(N,F) \, .
$$

\subsection{ \ }\label{ssec1.2}
To prove the first inequality in Theorem 1 we use a result of Borek \cite{Bo} which compares the successive minima and the slopes of hermitian vector bundles over $S$. Namely, according to \cite{Bo}, Th.~1, if $\sigma_k (\bar E)$ is the $k$-th slope of $\bar E$, the following inequality holds :
\begin{eqnarray}
\label{eq1}
&&0 \leq \mu_k (\bar E) + \sigma_k (\bar E)  \, . \nonumber
\end{eqnarray}
Similarly
\begin{eqnarray}
\label{eq2}
&&0 \leq \mu_{N+1-k} (\bar E^*) + \sigma_{N+1-k} (\bar E^*)  \, . \nonumber
\end{eqnarray}
On the other hand, we know that
$$
\sigma_k (\bar E) + \sigma_{N+1-k} (\bar E^*) = 0
$$
(see \cite{Gaudron}, 5.15(2)). So, by adding up, we get
\begin{eqnarray}
&&0 \leq \mu_k (\bar E) + \mu_{N+1-k} (\bar E^*) \, . \nonumber
\end{eqnarray}
\subsection{ \ }\label{ssec1.3}
The second inequality in Theorem~1 will be proved by reducing it to the case $F = {\mathbb Q}$. For every positive integer $k \leq Nr$ let $\lambda_k $ be the infimum of the set of real numbers $\lambda$ such that there exist $k$ vectors $e_1 , \ldots , e_k \in E$ which are ${\mathbb Q}$-linearly independent in $E \underset{\mathbb Z}{\otimes} {\mathbb Q}$ and such that, for every $\sigma \in \Sigma$ and every $i = 1,\ldots , k$,
$$
\Vert e_i \Vert_{\sigma} \leq \exp (\lambda) \, .
$$
The following lemma is used in \cite{Thunder}.

\bigskip

\noindent {\bf Lemma 1.} {\it For every positive integer $k \leq N$, the following inequality holds :}
$$
\mu_{k+1} (\bar E) \leq \lambda_{kr+1}  \, .
$$

\bigskip

\noindent {\it Proof.} Let $e_1 , \ldots , e_{kr+1} \in E$ be vectors which are ${\mathbb Q}$-linearly independent, and $V$ (resp. $W$) the $F$-vector space (resp. the ${\mathbb Q}$-vector space) spanned by these vectors. Since $W \subset V$ and $\dim_{\mathbb Q} (V) = r \dim_F (V)$ we get
$$
r \dim_F (V) \geq kr+1 \, ,
$$
hence $\dim_F (V) \geq k+1$.
\noindent The lemma follows from this inequality and the definition of successive minima.

\subsection{ \ }\label{ssec1.4}

Let 
$
E^{\vee} = {\rm Hom}_{\mathbb Z} (E , {\mathbb Z})
$
and
$
\omega = {\rm Hom}_{\mathbb Z} ({\mathcal O}_F , {\mathbb Z}) \, .
$
The morphism
$$
\alpha : E^* \otimes_{{\mathcal O}_F} \omega \to E^{\vee}
$$
mapping $u \otimes T$ to $u \circ T$ is an isomorphism of ${\mathcal O}_F$-modules. If ${\rm Tr} \in \omega$ is the trace morphism, we endow $\omega$ with the hermitian metric such that $\vert {\rm Tr} \vert_{\sigma} = 1$ (resp. $\vert {\rm Tr} \vert_{\sigma} = 2$) if $\sigma = \bar\sigma$ (resp. $\sigma \ne \bar\sigma$). For every $\sigma \in \Sigma$, the morphism
$$
E_{\sigma}^{\vee} \to E_{\sigma}^*
$$
induced by $\alpha$ is an isometry (\cite{Lattice}, p.~354). For any positive integer $k \leq Nr$, let $\lambda_k^{\vee}$ be the infimum of the set of real numbers $\lambda$ such that there exist $k$ vectors $e_1 , \ldots , e_k \in E^{\vee}$ which are linearly independent over ${\mathbb Q}$ and such that, for every $i = 1,\ldots , k$,
$$
\sum_{\sigma \in \Sigma} \Vert e_i \Vert_{\sigma} \leq \exp (\lambda) \, .
$$
According to \cite{Ban} Theorem~2.1 and section 3, we have, for $k=1,\ldots , Nr$,
\begin{equation}
\lambda_k + \lambda_{Nr+1-k}^{\vee} \leq \frac{3}{2} \log (Nr) \, .
\end{equation}

\subsection{ \ }\label{ssec1.5}

Since $\omega$ is invertible we have
$$
E^* \simeq E^{\vee} \otimes \omega^{-1}
$$
and, for any $v \in \omega^{-1}$, $v \ne 0$,
\begin{equation}
\mu_k (\bar E^*) \leq \mu_k (\bar E^{\vee}) + \underset{\sigma \in \Sigma}{\rm sup} \, \log \Vert v \Vert_{\sigma} \, .
\end{equation}
By Minkowski theorem we can choose $v$ such that, for every $\sigma \in \Sigma$,
$$
r \log \Vert v \Vert_{\sigma} \leq r \log (2) + \log {\rm covol} (\omega^{-1}) - \log {\rm vol} (B) \, ,
$$
where ${\rm vol} (B)$ is the volume of the unit ball in the real vector space $\omega_{\mathbb R}^{-1}$ and ${\rm covol} (\omega^{-1})$ is the covolume of the lattice $\omega^{-1}$. We have
$$
{\rm vol} (B) = 2^{r_1} \, \pi^{r_2}
$$
and, according to \cite{Lattice} p.~355,
$$
\log {\rm covol} (\omega^{-1}) = \log \vert \Delta_F \vert - 2r_2 \log (2) \, .
$$
So we can choose $v \in \omega^{-1}$, $v \ne 0$, such that
\begin{equation}
\underset{\sigma \in \Sigma}{\rm sup} \log \Vert v \Vert_{\sigma} \leq \frac{1}{r} \log \vert \Delta_F \vert - \frac{r_2}{r} \log (\pi) \, .
\end{equation}

\subsection{ \ }\label{ssec1.6}

>From Lemma~1 and the fact that
$$
\sum_{\sigma \in \Sigma} \Vert x \Vert_{\sigma} \leq r \, \underset{\sigma}{\rm sup} \, \Vert x \Vert_{\sigma}
$$
we get, for every $k \leq N$,
\begin{equation}
\mu_{k+1} (\bar E^{\vee}) \leq \lambda_{kr+1}^{\vee} + \log (r) \, .
\end{equation}
Therefore, using (3) and (4), we get
\begin{eqnarray}
&&\mu_k (\bar E) + \mu_{N+1-k} (\bar E^*) \nonumber \\
&\leq &\lambda_{(k-1)r+1} + \mu_{N+1-k} (\bar E^{\vee}) + \frac{1}{r} \log \vert \Delta_F \vert - \frac{r_2}{r} \log (\pi) \nonumber \\
&\leq &\lambda_{k+1-r} + \lambda_{(N-k)r+1}^{\vee} + \log (r) + \frac{1}{r} \log \vert \Delta_F \vert - \frac{r_2}{r} \log (\pi) \, . \nonumber
\end{eqnarray}
Since, by (2),
$$
\lambda_{k+1-r} + \lambda_{(N-k)r+1}^{\vee} \leq \lambda_{kr} + \lambda_{Nr-kr+1}^{\vee} \leq \frac{3}{2} \log (Nr) \, ,
$$
Theorem~1 follows.

\section{Lower bounds for the height of irreducible divisors}\label{sec2}

\subsection{ \ }\label{ssec2.1}
\hfill
Let $S = {\rm Spec} ({\mathcal O}_F)$ be as above. Consider a semi-stable curve $X$ over $S$ such that $X$ is regular and its generic fiber $X_F$ is geometrically irreducible of genus $g$. Let $h_X$ be an hermitian metric, invariant under complex conjugation, on the variety $X({\mathbb C})$ of complex points of $X$. Let $\omega_0$ be the associated K\"ahler form, defined by the formula
$$
\omega_0 = \frac{i}{2\pi} \, h_X \left( \frac{\partial}{\partial z} , \frac{\partial}{\partial z} \right) dz \, d\bar z
$$
if $z$ is any local holomorphic coordinate on $X({\mathbb C})$. Let $\bar L = (L,h)$ be an hermitian line bundle over $X$ (with $h$ invariant under complex conjugation). If $L_{\mathbb C}$ is the restriction of $L$ to $X({\mathbb C})$, the vector space $H^0 (X({\mathbb C}) , L_{\mathbb C})$ of holomorphic sections of $L_{\mathbb C}$ on $X({\mathbb C})$ is equipped with the sup norm
$$
\Vert s \Vert_{\rm sup} = \underset{x \in X({\mathbb C})}{\rm sup} \Vert s(x) \Vert \, ,
$$
where $\Vert \cdot \Vert$ is the norm defined by $h$, and with the $L^2$-norm
$$
\Vert s \Vert_{L^2}^2 = \sup_{\sigma} \int_{X_{\sigma}} \Vert s(x) \Vert^2 \, \omega_0 \, ,
$$
where $\sigma$ runs over all complex embeddings of $F$ and $X_{\sigma} = X \underset{{\mathcal O}_F}{\otimes} {\mathbb C}$ is the corresponding complex variety. We let
$$
A(\bar L_{\mathbb C}) = \sup_s \log (\Vert s \Vert_{\rm sup} / \Vert s \Vert_{L^2}) \, ,
$$
where $s$ runs over all sections of $L_{\mathbb C}$.

\smallskip

Consider the relative dualizing sheaf $\bar\omega_{X/S}$ of $X$ over $S$, equipped with the metric dual to $h_X$, and let $\bar M = \bar L \otimes \bar\omega_{X/S}^*$.
We endow the ${\mathcal O}_F$-module
$$
H^1 = H^1 (X , M^{-1})
$$
with the $L^2$-metric and we denote by $\mu_k (H^1)$ its successive minima, $k=1 , \ldots , N = \dim_F H^1 (X_F , M^{-1})$.

\smallskip

Let now $D$ be an irreducible divisor on $X$, flat over $S$, of degree $d$ on $X_F$. We are interested in the Faltings height $h_{\bar L} (D)$ of $D$ with respect to $\bar L$. Recall \cite{BGS} that $h_{\bar L} (D) \in {\mathbb R}$ is the arithmetic degree of the restriction of $\bar L$ to $D$.
Let $t = \dim_F H^0 (X_F , L(-D))$ and assume that $N > t$.

\bigskip

\noindent {\bf Proposition 1.} {\it The following inequality holds :}
$$
\frac{h_{\bar L} (D)}{dr} \geq \mu_{N-t} (H^1) - A (\bar L_{\mathbb C}) - C(N,F) \, .
$$

\bigskip

\noindent {\it Proof.} 
 To prove Proposition 1, let $s \in H^0 (X,L)$ be a section of $L$ which does not belong to the vector space $H^0 (X_F , L(-D))$. The restriction of $s$ to $D({\mathbb C})$ does not vanish hence, since $D$ is irreducible, for any point $P$ in $D({\mathbb C})$ we have $s(P) \ne 0$. The height of $D$ can be computed using $s$ (\cite{BGS} (3.2.2))
$$
h_{\bar L} (D) = h_{\bar L} ({\rm div} (s_{\vert D})) - \sum_{\alpha} \log \Vert s (P_{\alpha})\Vert \geq - \sum_{\alpha} \log \Vert s(P_{\alpha})\Vert \, ,
$$
where $D({\mathbb C}) = \underset{\alpha}{\sum} \, P_{\alpha}$. Next we have
$$
\sum_{\alpha} \log \Vert s(P_{\alpha}) \Vert \leq d r \log \Vert s \Vert_{\rm sup} \leq d r(\log \Vert s \Vert_{L^2} + A (\bar L_{\mathbb C})) \, .
$$

Let $\bar E = (H^0 (X,L) , h_{L^2})$. If $t$ is the rank of $H^0 (X,L(-D))$ we can choose $s$ such that
\begin{equation}
\log \Vert s \Vert_{L^2} \leq \mu_{t+1} (\bar E) \, .
\end{equation}
By Theorem~1
\begin{equation}
\mu_{t+1} (\bar E) \leq -\mu_{N-t} (\bar E^*) + C(N,F) \, ,
\end{equation}
and, by Serre duality, $\bar E^* = H^1 (X,M^{-1})$ with the $L^2$-metric. Therefore Proposition~1 follows from (6) and (7).

\subsection{ \ }\label{ssec2.2}

We keep the hypotheses of Proposition ~1 and we denote by $\bar M \cdot \bar M \in {\mathbb R}$ the self-intersection of the first arithmetic Chern class $\hat c_1 (\bar M) \in \widehat{\rm CH}^1 (X)$. Let $\delta = \deg (L)$ be the degree of $L$ on $X_F$ and $m = \deg (M) = \delta - 2g + 2$.

\bigskip

\noindent {\bf Proposition 2.} {\it Assume that $\delta$ is even and that
$$
2g+1 \leq d \leq \delta \leq 2d-2 \, .
$$
Then}
$$
\frac{h_{\bar L} (D)}{dr} \geq \frac{\bar M \cdot \bar M}{2mr} - A (\bar L_{\mathbb C}) - C (N,F) - \log (\delta (\delta - g + 1)) \, .
$$

\bigskip

\noindent {\it Proof.} According to \cite{Stable} Th.~2 and \cite{Stable} 2.3.1, the inequality
\begin{equation}
\mu_k (\bar E^*) \geq \frac{\bar M \cdot \bar M}{2mr} - \log (\delta (\delta - g + 1))
\end{equation}
holds 
$$
k \geq \frac{m}{2} + g = \frac{\delta}{2} + 1 \, .
$$
Consider the exact sequence of cohomology groups
\begin{eqnarray}
\label{eq5}
&&0 \to H^0 (X_F , L(-D)) \to H^0 (X_F , L) \to H^0 (D_F , L_{\vert D}) \nonumber \\
&&\to H^1 (X_F , L(-D)) \to H^1 (X_F , L) \, . 
\end{eqnarray}

\smallskip

We first assume that $\delta > d + 2g - 2$ i.e.
$$
\deg (L (-D)) > 2g-2 \, .
$$
This implies $H^1 (X_F , L(-D)) = 0$ and
$$
N-t = \dim_F H^0 (D_F , L_{\vert D}) = d \, .
$$
Since $d \geq \frac{\delta}{2} + 1$, the proposition follows from Proposition~1 and (8).

\smallskip

Next, we assume that
$$
d \leq \delta \leq d+2g-2 \, ,
$$
and we apply Clifford's theorem to the Serre dual of $L(-D)$ on $X_F$. It is special unless $H^0 (X_F , L(-D)) = 0$, in which case $t = 0$ hence
$$
N-t =\delta - g + 1 \geq \frac{\delta}{2} + 1
$$
since $\delta \geq 2g$, and we can conclude as above. 

\smallskip

When $H^0 (X_F , L(-D))$ does not vanish, Clifford's theorem says that
$$
\dim_F \, H^1 (X_F , L(-D)) - 1 \leq \frac{1}{2} \deg (\omega_{X/S} \otimes L^{-1} (D)) = g-1-\frac{\delta}{2} + \frac{d}{2} \, .
$$
>From (\ref{eq5}) it follows that
$$
N-t \geq d - \dim H^1 (X_F , L(-D))
$$
and therefore
$$
N-t \geq \frac{d}{2} + \frac{\delta}{2} - g \, .
$$
\smallskip

Since $d \geq 2g + 1$ this implies
$$
N-t \geq \frac{\delta}{2} + \frac{1}{2}
$$
and, since $\delta$ is even, we get
$$
N-t \geq \frac{\delta}{2} + 1 
$$
and the proposition follows from Proposition~1 and (8).

\subsection{ \ }\label{ssec2.3}

For any hermitian line bundle $\bar L$ on $X$, and any integer $d$, we define
$$
e(\bar L , d) = \inf_{\deg (D) = d} \frac{h_{\bar L} (D)}{d}
$$
and
$$
e(\bar L , \infty) = \lim_d \inf e(\bar L , d) \, .
$$

\bigskip

\noindent {\bf Theorem 2.} {\it If $\deg (L)$ is positive we have :}
$$
e(\bar L , \infty) \geq \frac{\bar L \cdot \bar L}{2 \deg (L)}  \, .
$$

\bigskip

\noindent {\it Proof.} By definition
$$
e(\bar L , \infty) = \lim_{n \to \infty} \inf_{\deg (D) = d \geq  n}\frac{h_{\bar L} (D)}{d} \,  .$$
Assume that $n \geq 2g + 1$ and $ n \geq \deg (L) + 3$ . Then, for any $d \geq n$,
there exists an even integer $k$ such that, if $\delta = k \deg (L)$, the inequalities
$$
2g+1 \leq d \leq \delta \leq 2d - 2 \, 
$$
hold. Fix a K\"ahler metric $h_X$ on $X({\mathbb C})$ (invariant by complex conjugation) and let 
$$
\bar M = \bar L^{\otimes k} \otimes \bar\omega^* \, .
$$
From Proposition~2 applied to $\bar L^{\otimes k}$ we get, for any irreducible horizontal divisor $D$ of degree $d$,
\begin{equation}
k \frac{h_{\bar L} (D)}{d \, r} \geq \frac{\bar M \cdot \bar M}{2 \deg (M) \, r} - A (\bar L_{\mathbb C}^{\otimes k}) - C (N , F) - \log (\delta (\delta - g + 1)) \, .
\end{equation}
When $n$ tends to infinity, the same is true for $d$ and $k$.  Therefore
\begin{equation}
\lim_{n \to \infty} \frac{\log (\delta (\delta - g + 1))}{k} = 0 \, .
\end{equation}
The rank $N$ of $H^{0} (X_F , L^{\otimes k})$ is $\delta - g + 1$ so, by (1), we have
\begin{equation}
\lim_{n \to \infty} \frac{C(N , F)}{k} = 0 \, . 
\end{equation}
According to a result of Gromov (\cite{ARR} Lemma 30) the quantity $\exp A (\bar L_{\mathbb C}^{\otimes k})$ is bounded from above by a polynomial in $k$. Therefore
\begin{equation}
\lim_{n \to \infty} \frac{A (\bar L_{\mathbb C}^{\otimes k})}{k} = 0 \, .
\end{equation}
Finally
$$
\deg (M) = k \deg (L) - 2g + 2
$$
and
$$
\bar M \cdot \bar M = (k \, \bar L - \bar\omega)^2 \, ,
$$
therefore
\begin{equation}
\lim_{n \to \infty} \frac{\bar M \cdot \bar M}{k \deg (M)} = \frac{\bar L \cdot \bar L}{\deg (L)} \, .
\end{equation}
The theorem follows from (10)--(14).

\subsection{ \ }\label{ssec2.4}

In \cite{Z} S.~Zhang defines
$$
e_{\bar L} = \underset{D}{\rm inf} \, \frac{h_{\bar L} (D)}{r \deg (D)}
$$
and
$$
e'_{\bar L} = \lim_D \inf \frac{h_{\bar L} (D)}{r \deg (D)} \, ,
$$
where $D$ runs over all irreducible horizontal divisors on $X$.

\bigskip

\noindent {\bf Lemma 2.} {\it When $\deg (L)$ is positive we have}
$$
e(\bar L,\infty) = r \, e'_{\bar L} \, .
$$

\bigskip

\noindent {\it Proof.} By definition
\begin{equation}
e(\bar L , \infty) = \lim_{n} \inf_{\deg (D) \geq n} \frac{h_{\bar L} (D)}{\deg (D)} \, .
\end{equation}
For any positive integer $n$ let $X(n)$ be the set of horizontal irreducible divisors $D$ such that
$$
\deg (D) < n \quad \mbox{and} \quad h_{\bar L} (D) \leq (e (\bar L , \infty) + 1) \, n \, .
$$
>From \cite{BGS}, Cor.~3.2.5, we know that $X(n)$ is finite and we get
\begin{equation}
r \, e' (\bar L) = \lim_{n} \inf_{D \notin X(n)} \frac{h_{\bar L} (D)}{\deg (D)} \, .
\end{equation}
The complement of $X(n)$ consists of those $D$ such that either $\deg (D) \geq n$ or $\deg (D) \leq n$ and $h_{\bar L} (D) > (e(\bar L , \infty) + 1) n$. In the second case we have
$$
\frac{h_{\bar L} (D)}{\deg (D)} > e (\bar L , \infty) + 1 \, .
$$
Therefore (16) and (17) imply
$$
r \, e' (\bar L) = {\rm Inf} (e(\bar L , \infty) , e (\bar L , \infty) + 1) = e(\bar L , \infty) \, .
$$
\hfill q.e.d.

\bigskip

When the first Chern form of $\bar L_{\mathbb C}$ is semi-positive and $\deg (L_{\vert E}) \geq 0$ for any vertical irreducible divisor $E$ on $X$, Theorem~6.3 in \cite{Z} states that
$$
r \, e'_{\bar L} \geq \frac{\bar L \cdot \bar L}{2 \deg (L)} \, .
$$
Therefore Theorem~2 is not new in that case.

\subsection{ \ }\label{ssec2.5}

We come back to the situation of \S~2.1 and 2.2, and we fix an integer $k \geq 1$. Furthermore we assume that the first Chern form of $\bar M_{\mathbb C}$ is positive and that $\deg (M_{\vert E}) \geq 0$ for any vertical irreducible divisor $E$ on $X$. If $k > 1$ define
$$
D(m,k) = (m+g) \sum_{\alpha = 0}^{{\rm Inf} (k-1 , g)} \begin{pmatrix} m+g-k-\alpha \\k-1-\alpha \end{pmatrix} \begin{pmatrix} g \\ \alpha \end{pmatrix} \, ,
$$
and let $D(m,1)=1$.

\bigskip

\noindent {\bf Proposition 3.} {\it Assume that $\delta \geq d \geq k$ and that either $m > 2k > 2$ or $m > k = 1$.
Then the following inequality holds :}
$$
\frac{h_{\bar L} (D)}{dr} \geq \frac{k}{m^2 \, r} \, \bar M^2 - \frac{2k}{m} \, e_{\bar M} +  \, e_{\bar M} - A (\bar L_{\mathbb C}) - C (N,F) - \frac{\log D(m,k)}{m^2} - 1 \, .
$$

\bigskip

\noindent {\it Proof.} According to \cite{Secant} Th.~4 i) (resp. \cite{Vanishing} Th.~2) we have
\begin{equation}
1 + \mu_k (H^1) \geq \frac{k}{m^2 \, r} \, \bar M \cdot \bar M - \frac{2k}{m} \, e_{\bar M} +  \, e_{\bar M} - \frac{\log D(m,k)}{m^2}
\end{equation}
as soon as $m > 2k > 2$ \footnote{Theorem~4, i) in \cite{Secant} assumes that $g \geq 2$ and the metric on $L_{\mathbb C}$ is admissible in the sense of Arakelov \cite{A}, but these extra hypotheses are not used in the proof of that statement.}(resp. $k=1$ and $m > 1$). If we assume that $\delta > d + 2g - 2$ we have $H^1 (X_F , L(-D)) = 0$ hence $N-t = d \geq k$. Therefore
$$
\mu_{N-t} (H^1) \geq \mu_k (H^1)
$$
and the proposition follows from (18) and Proposition~1. When $d \leq \delta \leq d + 2g - 2$ we consider the Serre dual of $L(-D)$ over $X_F$. It is special unless $t=0$, in which case
$$
N-t = \delta - g + 1 = m+g-1 \geq k \, .
$$
When $t \ne 0$, Clifford's theorem says that
$$
\dim H^1 (X_F , L(-D)) - 1 \leq \frac{1}{2} \deg (\omega \otimes L^{-1} (D)) = g - 1 - \frac{\delta}{2} + \frac{d}{2} \, ,
$$
and
$$
N-t \geq \frac{\delta}{2} + \frac{d}{2} - g \, .
$$
But
$$
\frac{\delta}{2} - g = \frac{m}{2} - 1 \geq k-1 \, ,
$$
hence
$$
N-t \geq k + \frac{d}{2} - 1
$$
and $N-t \geq k$ since $d \geq 1$.

\smallskip

Again, the proposition follows from (18) and Proposition~1.

\subsection{ \ }\label{ssec2.6}

We now assume that $g \geq 2$ and we let $\bar\omega$ be the relative dualizing sheaf $\omega_{X/S}$ of $X$ over $S$, equipped with its Arakelov metric \cite{A}. As in \ref{ssec2.3} above we consider
\begin{equation}
e(\bar\omega , d) = \underset{\deg (D)=d}{\rm inf} \frac{h_{\bar\omega} (D)}{d} \, .
\end{equation}

\bigskip

\noindent {\bf Theorem 3.} {\it There is a constant $C = C(g,r)$ such that the following inequalities hold:}
\begin{equation}
e(\bar\omega ,d) \geq \frac{\bar\omega \cdot \bar\omega}{4g(g-1)} \, \frac{dg + g-1}{d+2g-2} - \frac{g-1}{d+2g-2} \log \vert \Delta_F \vert - C \, \frac{\log (d)}{d} \, , 
\end{equation}
{\it and, if $d \geq 2g+1$,}
\begin{equation}
e(\bar\omega , d) \geq \frac{\bar\omega \cdot \bar\omega}{4(g-1)} \, \frac{d - 2g + 1}{d-g} - \frac{g-1}{d-g} \log \vert \Delta_F \vert - C \, \frac{\log (d)}{d} \, . 
\end{equation}

\bigskip

\noindent {\it Proof.} To prove (19) we apply Proposition~3 to a power $\bar L = \bar\omega^{\otimes n}$ of $\bar\omega$. We take $k=d$. When $d=1$, (19) follows from the inequalities
$$
e(\bar\omega , 1) \geq r \, e_{\bar\omega}
$$
and
\begin{equation}
r \, e_{\bar\omega} \geq \frac{\bar\omega \cdot \bar\omega}{4g(g-1)}
\end{equation}
(cf. \cite{E}). When $d>1$, the condition $m>2k$ in Proposition~3 becomes
$$
(n-1)(g-1) > d \, ,
$$
i.e.
$$
n > \frac{d}{g-1} + 1 \, .
$$
We take 
$$
n = \left[ \frac{d}{g-1} \right] + 2 \, .
$$
According to Proposition~3, for any irreducible horizontal divisor $D$ of degree $d$,
\begin{eqnarray}
\frac{h_{\bar L} (D)}{d} &\geq &k \, \frac{\bar\omega \cdot \bar\omega}{4 (g-1)^2} + r \, e_{\bar\omega} \left( n-1-\frac{k}{g-1} \right) \nonumber \\
&- &r \left( A(\bar L_{\mathbb C}) + C(N,F) + \frac{\log D(m,k)}{m^2} + 1 \right) \, . \nonumber
\end{eqnarray}
Using the lower bound (21) for $e_{\bar\omega}$ and the fact that
$$
h_{\bar L} (D) = n \, h_{\bar\omega} (D)
$$
we get
\begin{eqnarray}
e(\bar\omega , d) &\geq &\frac{\bar\omega \cdot \bar\omega}{4g(g-1)} \, \frac{k+n-1}{n} \nonumber \\
&- &\frac{r}{n} \left( A (\bar L_{\mathbb C}) + C(N,F) + \frac{\log D(m,k)}{m^2} + 1 \right) \, . 
\end{eqnarray}
Since
$$
n \leq 2 + \frac{d}{g-1}
$$
we get
\begin{equation}
\frac{k+n-1}{n} \geq \frac{dg + g - 1}{d + 2g - 2} \, .
\end{equation}
Gromov's estimate for $A (\bar\omega^{\otimes n})$ implies
\begin{equation}
\frac{A (\omega^{\otimes n})}{n} = O \left( \frac{\log (n)}{n} \right) = O \left( \frac{\log (d)}{d} \right) \, .
\end{equation}
From (1) we deduce that
\begin{equation}
\frac{r}{n} \, C(N,F) = \frac{1}{n} \log \vert \Delta_F \vert + O \left( \frac{\log (n)}{n} \right) \, .
\end{equation}
Finally, according to \cite{Secant} \S~3.8,
\begin{equation}
\log D(m,k) = O (m \log (m)) = O (d \log (d)) \, .
\end{equation}
The inequality (19) follows from (22)--(26).

\smallskip

To prove (20) we apply Proposition~2 to a power $\bar L = \bar\omega^{\otimes n}$ of $\bar\omega$. We get
\begin{equation}
e(\bar\omega , d) \geq \frac{n-1}{n} \, \frac{\bar\omega \cdot \bar\omega}{4(g-1)} - \frac{r}{n} (A (\bar L_{\mathbb C}) + C(N,F) + \log (\delta (\delta - g + 1)))
\end{equation}
as soon as
$$
2g+1 \leq d \leq (2g - 2) \, n \leq 2d-2 \, .
$$
We choose
$$
n = \left[ \frac{d-1}{g-1} \right] \geq \frac{d-g}{g-1}
$$
in which case
$$
\frac{n-1}{n} \geq \frac{d-2g+1}{d-g} \, .
$$
The second summand of the right-hand side of (27) is estimated as above. This proves (20).

\section{Upper bounds for the height of irreducible divisors}\label{sec3}

\subsection{ \ }\label{ssec3.1}

Let $X$ and $h_X$ be as in \S~2.1. Let $\bar L$ and $\bar M$ be two hermitian line bundles on $X$. We assume that $\deg (L) > 0$ and $\deg (L_{\vert E}) \geq 0$ for every vertical irreducible divisor $E$ on $X$. Let $D_0$ be an irreducible horizontal divisor,
$$
N = \dim_F H^0 (X_F , M)
$$
and
$$
t = \dim_F H^0 (X_F , M(-D_0)) \, .
$$
We assume that $N > t$. Denote by $\mu_k (H^1)$, $k = 1 , \ldots , N$, the successive minima of $H^1 = H^1 (X , \omega_{X/S} \otimes M^{-1})$ equipped with its $L^2$-metric. We write $\bar L \cdot \bar M \in {\mathbb R}$ for the arithmetic intersection of $\hat c_1 (\bar L)$ with $\hat c_1 (\bar M)$, and we write $D \pitchfork D_0$ to mean that $D$ is an irreducible horizontal divisor meeting $D_0$ properly.

\bigskip

\noindent {\bf Proposition 4.} {\it The following inequality holds :}
\begin{eqnarray}
\underset{D \pitchfork D_0}{\rm inf} \, \frac{h_{\bar L} (D)}{r \deg (D)} &\leq &\frac{\bar L \cdot \bar M}{r \deg (M)} - \mu_{N-t} (H^1) \frac{\deg (L)}{\deg (M)} \nonumber \\
&+ &\frac{\deg (L)}{\deg (M)} \, (A (\bar M_{\mathbb C}) + C (N,F)) \, . \nonumber
\end{eqnarray}

\bigskip

\noindent {\it Proof.} Let $\bar E = (H^0 (X,M) , h_{L^2})$ and choose a section $s \in H^0 (X,M)$ such that $s \notin H^0 (X_F,M(-D_0))$ and
$$
\log \Vert s \Vert_{L^2} \leq \mu_{t+1} (\bar E) \, .
$$
If ${\rm div} (s)$ is the divisor of $s$ we get (\cite{BGS} (3.2.2))
\begin{eqnarray}
\bar L \cdot \bar M &= &h_{\bar L} ({\rm div} (s)) - \int_{X({\mathbb C})} \log \Vert s \Vert \, c_1 (\bar L_{\mathbb C}) \nonumber \\
&\geq &h_{\bar L} ({\rm div} (s)) - r \deg (L) (\mu_{t+1} (\bar E) + A (\bar M_{\mathbb C})) \, . 
\end{eqnarray}
We can write
$$
{\rm div} (s) = \sum_{\alpha} D_{\alpha} + V
$$
where each $D_{\alpha}$ is irreducible and flat over $S$, and $V$ is effective and vertical on $X$. Therefore, by our assumption on $L$, we have
$$
h_{\bar L} ({\rm div} (s)) \geq \sum_{\alpha} h_{\bar L} (D_{\alpha})
$$
and
$$
\deg ({\rm div} (s)) = \sum_{\alpha} \deg (D_{\alpha}) \, .
$$
Therefore, since each $D_{\alpha}$ is transverse to $D_0$,
\begin{equation}
\frac{h_{\bar L} ({\rm div} (s))}{\deg (M)} \geq \underset{\alpha}{\rm inf} \, \frac{h_{\bar L} (D_{\alpha})}{\deg (D_{\alpha})} \geq \underset{D \pitchfork D_0}{\rm inf} \, \frac{h_{\bar L} (D)}{\deg (D)} \, .
\end{equation}
From Theorem~1 we get
\begin{equation}
\mu_{t+1} (\bar E) \leq -\mu_{N-t} (H^1) + C(N,F)
\end{equation}
and the proposition follows from (28), (29) and (30).

\subsection{ \ }\label{ssec3.2}

We keep the notation of the previous section and we let
$$
\bar K = \bar M \otimes \bar\omega_{X/S}^{*} \, , \ m = \deg (M) \quad \mbox{and} \quad d_0 = \deg (D_0) \, .
$$

\bigskip

\noindent {\bf Proposition 5.} {\it Assume that $m$ is even and
$$
2g+1 \leq d_0 \leq m \leq 2d_0 - 2 \, .
$$
The following inequality holds :}
\begin{eqnarray}
\underset{D \pitchfork D_0}{\rm inf} \, \frac{h_{\bar L} (D)}{r \deg (D)} &\leq &\frac{\bar L \cdot \bar M}{rm} - \frac{\bar K \cdot \bar K}{2r\deg (K)} \frac{\deg (L)}{m} \nonumber \\
&+ &\frac{\deg (L)}{m} (A(\bar M_{\mathbb C}) + C(N,F) + \log (m(m-g+1))) \, . \nonumber
\end{eqnarray}

\bigskip

\noindent {\it Proof.} The number $\mu_{N-t} (H^1)$ can be estimated from below using \cite{Stable} exactly as in the proof of Proposition~2. Therefore the proposition follows from Proposition~4.

\subsection{ \ }\label{ssec3.3}

Let $\bar L$ be an hermitian line bundle on $X$ such that $\deg (L) > 0$ and $\deg (L_{\vert E}) \geq 0$ for any irreducible vertical divisor $E$ on $X$. For any integer $d_0 \geq 1$ consider
$$
e' (\bar L , d_0) = \underset{D_0}{\rm sup} \underset{D \pitchfork D_0}{\rm inf} \, \frac{h_{\bar L} (D)}{\deg (D)} \, ,
$$
where $D_0$ runs over all irreducible horizontal divisors of degree $d_0$. Let
$$
e' (\bar L , \infty) = \lim_{d_0} \sup e' (\bar L , d_0) \, .
$$

\bigskip

\noindent {\bf Theorem 4.} {\it The following inequality holds :}
$$
e' (\bar L , \infty) \leq \frac{\bar L \cdot \bar L}{2 \deg (L)} \, .
$$

\bigskip

\noindent {\it Proof.} As in the proof of Theorem 2, when the integer $n$ is big enough,
for any $d_0 \geq n$ we can choose an even power $\bar M$ of $\bar L$
such that, if $m = \deg (M)$, the following inequalities hold :
$$
2g+1 \leq d_0 \leq m \leq 2d_0 - 2 \, .
$$
 Then we apply Proposition~5 to $\bar L$ and $\bar M$. If $\bar K = \bar M \otimes \bar\omega_{X/S}^{*}$ we get
\begin{equation}
\lim_{n \to \infty} \frac{\bar K \cdot \bar K}{\deg (K)} \, \frac{\deg (L)}{m} = \frac{\bar L \cdot \bar L}{\deg (L)}
\end{equation}
and
\begin{equation}
\lim_{nÊ\to \infty} \frac{\bar L \cdot \bar M}{m} = \frac{\bar L \cdot \bar L}{\deg (L)} \, .
\end{equation}
By the same estimates as in the proof of Theorem~2 we get
\begin{equation}
\lim_{n \to \infty} (A(\bar M_{{\mathbb C}}) + C(N,F) + \log (m (m - g + 1)))/m = 0 \, .
\end{equation}
The theorem follows from (31), (32), (33) and Proposition~5.

\bigskip

\noindent {\bf Remark.} For any $d_0$ we have   
$$
r \, e_{\bar L} \leq e' (\bar L , d_0)  \, .
$$
Therefore Theorem~3 implies
$$
r \, e_{\bar L} \leq \frac{\bar L \cdot \bar L}{2 \deg (L)} \, .
$$
But it does not follow from \cite{Z}, Th.~6.3.

\subsection{ \ }\label{ssec3.4}

We come back to the notation of \ref{ssec3.2} and we let
$$
k = \deg (K) = m - 2g + 2 \, .
$$
We fix an integer $h \geq 1$. We assume that the first Chern form of $\bar K_{\mathbb C}$ is positive and that $\deg (K_{\vert E}) \geq 0$ for every irreducible vertical divisor $E$ on $X$.

\bigskip

\noindent {\bf Proposition 6.} {\it Assume that $m \geq d_0 \geq h$ and that either $k > 2h > 2$ or $k > h = 1$. Then the following inequality  :}
\begin{eqnarray}
\underset{D \pitchfork D_0}{\rm inf} \, \frac{h_{\bar L} (D)}{r \deg (D)} &\leq &\frac{\bar L \cdot \bar M}{rm} - \frac{\deg (L)}{m} \left( \frac{h}{k^2r} {\bar K}^2-  \, \frac{2h}{k} e_{\bar K}+ e_{\bar K} \right)  \\
&+ &\frac{\deg (L)}{m} \left( A (\bar M_{\mathbb C}) + C(N,F) + \frac{\log D(k,h)}{h^2} + 1 \right) \, . \nonumber
\end{eqnarray}

\noindent {\it Proof.} This inequality follows from Proposition~4 by bounding $\mu_{N-t} (H^1)$ from below in the same way as in the proof of Proposition~3.

\subsection{ \ }\label{ssec3.5}

Assume now that $g \geq 2$ and let $\bar\omega$ be $\omega_{X/S}$ with its Arakelov metric. Recall that
$$
e' (\bar\omega , d_0) = \underset{\deg (D_0) = d_0}{\rm sup} \ \underset{D \pitchfork D_0}{\rm inf} \, \frac{h_{\bar L} (D)}{\deg (D)} \, .
$$

\bigskip

\noindent {\bf Theorem 5.} {\it There exists a constant $C = C(g,r)$ such that the following inequalities hold :}
\begin{eqnarray}
e' (\bar\omega , d_0) \leq \frac{\bar\omega \cdot \bar\omega}{4(g-1)} + \frac{2g-1}{4g(d_0 + 2g-2)} \, \bar\omega \cdot \bar\omega 
 + \, \frac{g-1}{d_0 + g-1} \log \vert \Delta_F \vert + C \, \frac{\log (d_0)}{d_0} \, , 
\end{eqnarray}
{\it and, when $d_0 \geq 2g+1$,}
\begin{eqnarray}
e' (\bar\omega , d_0) \leq \frac{\bar\omega \cdot \bar\omega}{4(g-1)} + \frac{\bar\omega \cdot \bar\omega}{4(d_0 - g)} + \frac{g-1}{d_0 - g} \log \vert \Delta_F \vert + C \, \frac{\log (d_0)}{d_0} \, . 
\end{eqnarray}

\bigskip

\noindent {\it Proof.} To prove (35) we apply Proposition~6 with $\bar L = \bar\omega$, $\bar M = \bar\omega^{\otimes n}$ and $h=d_0$.
When $d_0 = 1 < k$ we have $n (g-1) \geq g$.
When $d_0 > 1$ and
$$
k=n(2g-2) - 2g+2 > 2 \, d_0
$$
we get $n(g-1) > d_0 + g - 1$. 

\smallskip

In both cases we choose
$$
n=2 + \left[ \frac{d_0}{g-1} \right] \, .
$$
The right hand side of (34) (Proposition~6) becomes $X_1 + X_2$, with
$$
X_1 = \frac{n \, \bar\omega \cdot \bar\omega}{rn(2g-2)} - \frac{1}{n} \left( d_0 \, \frac{\bar\omega \cdot \bar\omega}{(2g-2)^2r} + \left( 1-\frac{2 \, d_0}{(n-1)(2g-2)} \right) (n-1) \, e_{\bar\omega} \right)
$$
and
$$
X_2 = \frac{\deg (L)}{m} \left( A(\bar M_{\mathbb C}) + C(N,F) + \frac{\log D (k,h)}{h^2} + 1 \right) \, .
$$
As in the proof of Theorem~3 we get
$$
X_2 \leq C \, \frac{\log (d_0)}{d_0} + \frac{1}{nr} \log \vert \Delta_F \vert
$$
and
$$
\frac{1}{n} \leq \frac{g-1}{d_0 + g - 1} \, .
$$
On the other hand, since 
$$r e_{\bar\omega} \geq \frac{\bar\omega \cdot \bar\omega}{4g(g-1)}\, ,$$
we get
\begin{eqnarray}
r \, X_1 &\leq &\bar\omega \cdot \bar\omega \left( \frac{1}{2g-2} - \frac{d_0}{n(2g-2)^2} - \frac{n-1}{4g(g-1)n} + \frac{d_0}{4ng (g-1)^2} \right) \nonumber \\
&= &\frac{\bar\omega \cdot \bar\omega}{4g(g-1)} \left( 2g-1 - \frac{d_0 - 1}{n} \right) \, . \nonumber
\end{eqnarray}
Since $n \leq 2 + \frac{d_0}{g-1}$ we get
\begin{eqnarray}
r \, X_1 &\leq &\frac{\bar\omega \cdot \bar\omega}{4g(g-1)} \left( 2g-1 - \frac{(d_0 - 1)(g-1)}{2g-2 + d_0} \right) \nonumber \\
&= &\frac{\bar\omega \cdot \bar\omega}{4(g-1)} + \frac{2g-1}{4g (d_0 + 2g-2)} \, \bar\omega \cdot \bar\omega \, . \nonumber
\end{eqnarray}
This proves (35).

\smallskip

To prove (36) we apply Proposition~5 when $\bar L = \bar\omega$ and $\bar M = \bar\omega^{\otimes n}$. If $d_0 \leq m \leq 2d_0 - 2$ we get 
$$
e(\bar L , d_0) \leq r Y_1 + r Y_2
$$ 
where
\begin{eqnarray}
Y_2 &= &\frac{\deg (L)}{m} \, (A(\bar M_{\mathbb C}) + C(N,F) + \log (m(m - g+1))) \nonumber \\
&\leq &C \, \frac{\log (d_0)}{d_0} + \frac{1}{nr} \log \vert \Delta_F \vert \nonumber
\end{eqnarray}
as in the proof of Theorem~3, and
\begin{eqnarray}
r \, Y_1 &= &\frac{\bar L \cdot \bar M}{m} - \frac{\bar K \cdot \bar K}{2 \deg (K)} \, \frac{\deg (L)}{m} \nonumber \\
&= &\frac{\bar\omega \cdot \bar\omega}{2g-2} - \frac{n-1}{4n(g-1)} \, \bar\omega \cdot \bar \omega \nonumber \\
&= &\frac{\bar\omega \cdot \bar\omega}{4(g-1)} + \frac{\bar\omega \cdot \bar\omega}{4n(g-1)} \, . \nonumber
\end{eqnarray}
Since $n(g-1) \leq d_0 - 1$ we can assume that
$$
n =\left[ \frac{d_0 - 1}{g-1} \right] \, ,
$$
hence $n \geq \frac{d_0 - 1}{g-1} - 1$. This implies
$$
\frac{1}{n} \log \vert \Delta_F \vert \leq \frac{g-1}{d_0 - g} \log \vert \Delta_F \vert
$$
and
$$
r \, Y_1 \leq \frac{\bar\omega \cdot \bar\omega}{4(g-1)} + \frac{\bar\omega \cdot \bar\omega}{4(d_0 - g)} \, ,
$$
from which (36) follows.

\end{document}